\documentclass[12pt]{article}
\usepackage{amsmath,amsthm, amssymb,latexsym,amscd}

\theoremstyle{plain}
\newtheorem{theorem}{Theorem}[section]
\newtheorem{proposition}[theorem]{Proposition}
\newtheorem{corollary}[theorem]{Corollary}
\newtheorem{lemma}{Lemma}[section]
\newtheorem{sublemma}{Sublemma}[section]

\newtheorem{definition}{Definition}[section]

\renewcommand{\mod}{{\mathop{\mbox{{\rm mod}}}}}

\newcommand{\ie}{ {\textit{ i.e.}} }

\DeclareMathOperator{\comp}{Comp}

\newenvironment{plm}{\noindent {\it Proof.\ }} { \hfill\qed\\ }
\newenvironment{prf}{\noindent {\it Proof.\ }} { \hfill\qed\\ }
\newenvironment{proofof}[1]{\medskip 
\noindent{\it Proof of #1.}}{ \hfill\qed\\ }

\newcommand{\LL}{\mathcal{L}}
\newcommand{\Comp}{\comp}

\newcommand{\A}{\mbox{\bf A}}

\newcommand\R{{\mathbb R}}
\newcommand\C{{\mathbb C}}
\newcommand\D{{\mathbb D}}

\newcommand{\Crit}{\mbox{Crit}}

\newcommand{\F}{\mbox{\bf F}}
\newcommand{\I}{\mbox{\bf I}}
\newcommand{\II}{{\small \mbox{\bf I}}}

\newcommand\N{{\mathbb N}}

\newcommand{\BP}{\mbox{\bf P}}
\newcommand{\Q}{\mathbb Q}
\newcommand{\Part}{\mathcal{P}}
\newcommand{\PP}{{\mathbb E}}
\newcommand{\QQ}{{\mathbb F}}

\newcommand{\U}{\mbox{\bf U}}

\newcommand{\Z}{\mbox{\bf Z}}

\newcommand\diam{\mathop{{\rm diam}}}

\newcommand{\dom}{\mathop{{\rm Dom}}}

\begin{document}
\title{Local connectivity and quasi-conformal rigidity of non-renormalizable  polynomials}
\author{Oleg Kozlovski and Sebastian van Strien}
\date{19th Sept 2006 (updated 2009)}
\maketitle

\begin{abstract} We prove that topologically conjugate non-renormalizable
polynomials are quasi-conformally conjugate. 
From this we derive that each such polynomial can be approximated
by a hyperbolic polynomial. 
As a by-product we prove that the Julia set of a
non-renormalizable polynomial
with only hyperbolic periodic points is locally connected,
and the Branner-Hubbard conjecture.
The main tools are the enhanced nest construction (developed
in a previous joint paper with Weixiao Shen \cite{KSS1})
and a lemma of Kahn and Lyubich (for which we give an elementary
proof in the real case). 
\end{abstract}

\section{Statement of Main Results}
Let $f\colon \C\to \C$ be a polynomial.
Its filled Julia set $K(f)$ is the set of points which do not escape to infinity, i.e. all $z$ so that $|f^n(z)|\not\to \infty$, and its Julia set
$J(f)$ is equal to $\partial K(f)$. We call such a polynomial $f$
{\em hyperbolic}, if each of its critical points is contained in the basin
of a hyperbolic periodic attractor.  A classical question is whether
one can approximate each polynomial by a hyperbolic polynomial
of the same degree (i.e. only small changes of the coefficients
are needed to 'make' the polynomial hyperbolic).
In this paper we give an answer to this question under the additional assumption that the original polynomial is non-renormalizable (or more, generally,  only finitely renormalizable) and has only hyperbolic periodic points. We derive this result
by proving that topologically conjugate non-renormalizable polynomials (without only hyperbolic periodic orbits) are
quasi-conformally conjugate,
while dealing with the
additional complication that there may be several critical points.
In addition we will show that under the same non-renormalizability assumption  each point in the Julia set is contained in arbitrarily small puzzle pieces (and therefore that the Julia set is locally connected) and that there is no measurable invariant linefield supported on the Julia set. 

Let us define these notions now. Since $f$ is a polynomial,
it is conformally conjugate near $\infty$ to $z\mapsto z^d$ where
$d$ is the degree of $f$. 
 curve
which under  this conformal equivalence near $\infty$ 
is mapped to a circle (respectively
to  a line-segment  through $0$ with rational angle) is called an 
{\em equipotential} (respectively a piece of a {\em periodic ray}). 
Iterates of an equipotential are all disjoint.
If $f$ has a connected Julia set, 
then the backward orbit of a piece of a periodic ray forms a curve $\gamma$ such that $f^k(\gamma)=\gamma$
and such that $\gamma$ lands on a periodic point; moreover
each periodic point has a periodic
ray landing on it, see for example \cite{MR2193309}.
If $f$ has also periodic attractors in $\C$ or $f$ has a Julia
set which is not connected, then we also 
associate to these orbits equipotentials and periodic rays 
as in Subsection~\ref{subsec:puzzles}.
A {\em puzzle piece} $V$ is a bounded set whose boundary consists
of pieces of equipotentials and periodic rays such that
when $z\in \partial V$ then $f^n(z)\notin \mbox{int}(V)$
for all $n\ge 0$. We say that a periodic point is {\em hyperbolic}
if its multiplier is not on the unit circle, and {\em repelling}
if its multiplier is outside the unit circle. 
For quadratic maps with only repelling periodic points
it is well-known how to construct puzzle pieces; 
in Subsection~\ref{subsec:puzzles} we extend this to general polynomials (with only hyperbolic periodic points). 

Given a puzzle piece $V$ containing a critical point 
$c$, define  $R_V$ be the first return map to 
$V$ and let $K_c(R_V)$ be the component of 
$K(R_V)=\{z; R^k(z)\in V\mbox{ for all }k\ge 0\}$
containing $c$.  
We say that $f$ is {\em finitely renormalizable} if there exists an integer $s$  such that for each critical point $c$ of $f$ 
and each puzzle piece $V\ni c$ for which the return time 
of $c$ to $V$ is  greater or equal to $s$, the set $K_c(R_V)$ is non-periodic.
We say that $f$ is {\em non-renormalizable}, if $s=2$. 

As usual, we often abbreviate the term 
`quasi-conformal conjugacy'  to `qc-conjugacy'.
Finally, we say that $f$ has an {\em invariant linefield} on 
$J(f)$, if there exists a measurable subset $E\subset J(f)$  and a measurable map which associates to Lebesgue
almost every $x\in E$ a line $l(x)$ through $x$
which is $f$-invariant in the sense that
$l(f(x))=Df(x)\,l(x)$. (So the absence of linefields
is obvious if the Julia set has zero Lebesgue measure.)

Our first result generalizes the celebrated result of 
\cite{Yoccoz} for finitely renormalizable quadratic maps,
and of \cite{AKLS} for finitely renormalizable unicritical maps.

\begin{theorem}[Local connectivity and absence of linefields]
\label{thm:main}
Assume that $f$ 
is a  finitely 
renormalizable polynomial whose periodic points 
are all hyperbolic. Then 
\begin{enumerate}
\item each point of its Julia set is contained in
arbitrarily small puzzle pieces;
\item if the Julia set is connected, then it is locally connected;
\item $f$ has no invariant linefields on its Julia set;
\end{enumerate}
\end{theorem}

Note that the puzzle construction we consider in the above theorem, 
also consists of rays and equipotentials
contained in basins of periodic attractors, see 
Subsection~\ref{subsec:puzzles}.



This theorem 
implies the following conjecture of Branner and Hubbard
(stated in the cubic case in \cite{MR1194004}).
(While this paper was in the final stages of
preparation we received a preprint 
W. Qiu and Y. Yin \cite{QiuYin} also proving this conjecture.)

\begin{corollary}[Branner - Hubbard conjecture]
The Julia set of a polynomial is totally disconnected if and only if
each critical component of the filled Julia set is aperiodic.
\end{corollary}
\begin{proof}
The only if part is obvious. 
If each critical component of the filled 
Julia set is aperiodic, then $f$ is finitely renormalizable
but where in the above definition 
we only need to consider puzzle pieces
bounded by preimages of equipotentials. 
But this implies that in the proof of the above theorem
all puzzle pieces are also of this type.
\end{proof}

\begin{theorem}[Rigidity of non-renormalizable polynomials]
\label{thm:rigidity}
Assume that $f$  is a  finitely 
renormalizable polynomial whose periodic points 
are all hyperbolic.
Moreover, assume that one of the following conditions hold
\begin{description}
\item[a)\quad] $f$ and $\tilde f$ are
topologically conjugate;
\item[b)\quad] $f,\tilde f$ are combinatorially equivalent 
and restricted  to their Fatou sets there exists a
compatible conjugacy 
between  $f$ and $\tilde f$
(these notions are explained in Definition~\ref{defcombequiv}) . 
\end{description}
Then $f$ and $\tilde f$ are quasi-conformally conjugate. 
\end{theorem}

For real polynomials one can drop the assumption that the map is
finitely renormalizable: this was shown for real quadratic maps in
\cite{GraSwi_hyper} and  \cite{Lyubich_dynI+II} and 
for general real polynomials  in \cite{KSS1}.
For finitely renormalizable complex quadratic maps, the above
result was shown by \cite{Yoccoz} and  for finitely renormalizable unicritical polynomials in \cite{AKLS}. 
Here combinatorial equivalence is defined as below (\cite{McM}); we 
note that \cite{LePrz} gives the construction of rays associated 
 to a super-attracting periodic point (in $\C$ or $\infty$), even if
 the Julia set is disconnected.

\begin{definition}[Combinatorially Equivalence
and Compatible Conjugacy]\label{defcombequiv} 
{\rm Let us say that $\alpha\sim_f \beta$ for 
rational numbers $\alpha,\beta$  
if the external rays of $f$ with angles $\alpha,\beta$
land on the same (pre-)periodic point. We then 
say that two polynomials $f,\tilde f$
are {\em combinatorially equivalent} if $\sim_f$ and $\sim_{\tilde f}$
define the same equivalence relation on the set of rational numbers.

A topological conjugacy $h$ between $f,\tilde f$ on their Fatou sets is called
{\em compatible} if for each periodic component $B$
of the Fatou set, $h$ extends continuously to the closure of $B$
and the correspondence between
periodic points in $\partial B$ and periodic points in $\partial h(B)$ 
coming from this extension of $h$, agrees with the one induced  
by the combinatorial equivalence.}
\end{definition}

By Lemma~\ref{lem:basins} below, a conjugacy on 
a periodic component $B$ of the Fatou set can always
be extended to $\partial B$  (for polynomials
which are non-renormalizable and with only hyperbolic
periodic points); it is not necessarily compatible
(since there are several homeomorphisms conjugating
the restriction of $z\mapsto z^d$ to $\{z\in \C; \, |z|=1\}$
to itself.) 
However,  if $f$ and $\tilde f$ have no escaping critical 
points and no periodic attractors, 
then the B\"ottcher coordinates at infinity
induce a (conformal) compatible conjugacy on the basin of infinity,
and so the compatibility in assumption (2) 
in Theorem~\ref{thm:rigidity} is automatic.

Also note that if $f,\tilde f$ are topologically conjugate and their
Julia sets are connected, they
are combinatorially equivalent:
the image of a periodic 
ray of $f$ by the conjugacy is a periodic curve which 
lands on a periodic point of $\tilde f$, and thereby 
by Theorem 2 of \cite{LePrz} homotopy 
equivalent relative to $\C\setminus K(\tilde f)$  to a periodic ray of $\tilde f$.
If the Julia set is non-connected, rays need not be smooth, as
they can split at critical points. Because of this the
connectedness assumption is needed for the implication
asserted in this paragraph:  all quadratic maps
with escaping critical points are topologically conjugate 
but they need not be combinatorially equivalent.

\bigskip

The following theorem was proved for finitely renormalizable 
quadratic maps by Yoccoz \cite{Yoccoz} and
for arbitrary real polynomials in \cite{KSS2}.

\begin{theorem}[Approximating finitely renormalizable by hyperbolic polynomials]\label{thm:density}
  Any finitely renormalizable polynomial $f$ 
  which has only  hyperbolic
  periodic points
  can be approximated by a hyperbolic polynomial $g$ of the same degree. If the Julia set of $f$ is connected, then $g$ can be 
  chosen so that its Julia set is also connected.
\end{theorem}

As mentioned before, 
we say that a polynomial is {\em hyperbolic} if all its
critical points are in basins of hyperbolic periodic attractors.
\medskip

The strategy of the proof is to use the enhanced nest construction
in \cite{KSS1} combined with a lemma due to Kahn and Lyubich
to obtain  complex bounds. As a by-product of the present paper, 
we obtain a significant simplification of the 
proof of the Key Lemma from \cite{KSS1}:
the  main technical part 
of \cite{KSS1} (Section 8-11)  can be replaced by this paper.

In order to prove the above theorems we shall need to work
with so-called complex box mappings. For unicritical polynomials the
construction of a complex box mapping is standard: It is a map $F: U\to V$,
where $V$ is a topological disk, $U$ is a union of (possibly
infinitely many) simply connected domains, all components of $U$
except one are mapped univalently onto $V$ by $F$ and $F$ restricted
on that other component of $U$ is a cover map onto $V$. When we
consider a multicritical polynomial we cannot construct a complex 
box mapping
which has the same nice properties as in the unicritical case. 
If the domain $V$
is connected and small, then it is possible that the first return map
to $V$ has infinitely many critical points, which of course 
is undesirable. So, the domain $V$ in our construction will have several connected
components. In the unicritical case all components of $U$ are
compactly contained in $V$, in the 
multicritical case we have to drop this
property as well.

\begin{definition}[Complex box mappings]\label{def1} 
  We say that a
  holomorphic map $F\colon U\to V$ 
  between open sets $U \subset V$ in
  $\C$ is a {\em complex box mapping} if the following hold:
  \begin{enumerate}
  \item  $F$ has finitely many critical points;
  \item $V$ is a  union of finitely many pairwise disjoint Jordan disks;
  \item every connected component $V'$ of $V$ is either a connected
    component of $U$ or the intersection of $V'$ and $U$ is a union of
    Jordan disks with pairwise disjoint closures which are compactly
    contained in $V'$,
  \item for each component $U'$ of $U$, $F(U')$ is
    a component of $V$.
  \end{enumerate}
\end{definition}

This  generalizes the well-known notion of a {\em 
(generalized) polynomial-like map}: in that case $V$ has only one component and
each component of $U$ is compactly contained in $V$.
Many authors also require that in addition $U$ has only a finite number of components.

Given a complex polynomial we will construct an induced complex box
mapping  $F : U \to V$ for this polynomial which has a few particular
properties (however, we will not need these properties in what follows).
Suppose that all critical points which are in the Julia set of the
polynomial are recurrent to each other. Then the domain $V$ can be
decomposed as $V=V_0 \cup V_1 \cup \cdots \cup V_k$, where domains
$V_1,\ldots, V_k$ are also components of $U$ and $U\setminus
(V_1\cup\cdots\cup V_k) \subset V_0$. Each component $V_i$,
$i=0,\ldots,k$, contains at least one critical point.

The filled Julia set of $F$ is then defined by
$K(F)=\{z\in U: f^n(z)\in U \mbox{ for any }n\in\N\}$
and  its {\it Julia set} as $J(F)=\partial K(F)$.
We denote by $Crit(F)$ the set of  critical points of $F$
and let $PC(F)=\bigcup_{n\ge 0 }f^n(Crit(F))$.
A {\em puzzle piece} is a component of $F^{-n}(V)$ for some
$n\geq 0$.  If it  contains a critical point then it is called a {\em critical puzzle piece}. 

Note that we allow there to be several critical points
in a component of the domain, and also that 
$F$ has no critical points.

\begin{definition}[Renormalizable complex box mapping]
The complex box mapping $F$ is called {\em renormalizable}
if there exists $s>1$ (called the {\em period}) and 
a puzzle piece $W$ 
containing a critical point $c$ of $f$ such that
$f^{ks}(c)\in W$, $\forall k\ge 0$ 
and such that $s$ is the first return time of $c$ to $W$.
\end{definition}

\newpage 

\begin{definition}[Itinerary of a puzzle piece relative to some curve family $\Gamma$]
Let $X$ be a collection of points in $\partial V$ intersecting each
component of $\partial V$,  so that for each point $y\in F^{-1}(X)$
there exists a simple curve in $V\setminus (U\cup PC(F))$ 
connecting $y$ to a point of the set $X$.
Let $\Gamma$ denote this collection 
of curves. For each component $U'$ of $F^{-n}(U)$  there exists a curve
connecting $\partial U'$ to $x$ of the form  $\gamma_{0}\cdots  \gamma_n$
where $F^k(\gamma_k)\in \Gamma$. The word $( \gamma_0,F(\gamma_1),\dots,F^n(\gamma_n))$
is called a {$\Gamma$-itinerary} of $U'$. 
A component of $F^{-n}(U)$ can have several
$\Gamma$-itineraries,  but different components of $F^{-n}(U)$ have different $\Gamma$-itineraries (as can be easily checked).
\end{definition}

We say that two complex box mappings are combinatorially equivalent,
if their critical points have the same itineraries:

\begin{definition}[Combinatorial equivalence of complex box mappings]
\label{def:equivbox}
Non-renormalizable complex box mappings 
$F\colon U\to V$ and $\tilde F\colon \tilde U\to \tilde V$ are called {\em combinatorially equivalent} w.r.t. some 
homeomorphism $H\colon V\to \tilde V$ with $H(U)=\tilde U$
and with $H(PC(F)\setminus U)=PC(\tilde F)\setminus \tilde U$
if there exists a curve family $\Gamma$ as in the previous definition, so that each 
critical point $c\in Crit(F)$ corresponds to a unique critical point 
$\tilde c\in \Crit(\tilde F)$
with the property that for every integer $k,n\ge 0$, the $\Gamma$-itineraries  
of the  component of $F^{-n}(U)$ containing $F^k(c)$ agree with the
$\tilde \Gamma$-itineraries
of the component of $\tilde F^{-n}(\tilde U)$ containing $\tilde F^k(\tilde c)$.
Here $\tilde \Gamma=H(\Gamma)$.
\end{definition}

For our purposes we will only need to work with non-renormalizable
complex box mappings. 
By defining the notion of renormalization of a complex box mapping,
it is not hard to extend the following theorem also to
finitely renormalizable complex box mappings.

\begin{theorem}[Rigidity for complex box mappings]\label{thm:main2}
Assume that $F\colon U\to V$  
is a  non-renormalizable complex box mapping
whose periodic points  are all repelling. Then 
\begin{enumerate}
\item each point of its Julia set is contained
in arbitrarily small puzzle pieces;
\item $F$ has no measurable invariant linefields on its Julia set;
\item Assume that $\tilde F\colon \tilde U\to \tilde V$ 
is another complex box mapping for which there exists
a quasi-conformal homeomorphism $H\colon \C\to \C$ so that
$H(V)=\tilde V$, $H(U)=\tilde U$, $\tilde F\circ H=H\circ F$ on $\partial U$
and so that $\tilde F$ is combinatorially equivalent to $F$ w.r.t. $H$. 
Moreover, assume that the boundary of each component of $U,V,\tilde U,\tilde V$
consists of piecewise smooth arcs. 
Then $F$ and $\tilde F$ are quasiconformally conjugate.
\end{enumerate}
\end{theorem}

The maps $F$ and $\tilde F$ can have several critical points and are not necessarily real.
Also note that it is not necessary to assume that $\tilde F$ has only hyperbolic periodic points.

The organisation of this paper is as follows.
In the first part of the paper
we shall show that rigidity of box mappings (Theorem~\ref{thm:main2}) implies the other
theorems. From Section~\ref{sec:rigboxbounds} we will prove
Theorem~\ref{thm:main2} using
the enhanced nest construction from \cite{KSS1} and a Lemma
by Kahn and Lyubich.

The authors would like to warmly thank Weixiao Shen.
At the end of his stay in Warwick,  in the spring of 2005,
Weixiao indicated that the Kahn-Lyubich lemma combined
with the enhanced nest construction from \cite{KSS1} gives complex bounds, see 
Section~\ref{sec:complexboundsenhanced}.    The authors also 
gratefully acknowledge several very useful discussions 
with Genadi Levin, in particular on the puzzle construction 
in Section~\ref{subsec:puzzles} and with Tan Lei on the definition of
combinatorial equivalence at the beginning of the paper.

\section{Local connectivity and absence of linefields (Theorem~\ref{thm:main2} implies 
Theorem~\ref{thm:main})}
\label{sec:polfrombox}
First we need to associate a puzzle partition
to any polynomial $f$ which only has hyperbolic periodic points,
and then use this to construct a complex box mapping
$F\colon U\to V$.
If $f$ has only repelling periodic points, then the construction is  
the multi-critical analogue of the usual Yoccoz puzzle partition, but  otherwise the construction will make use of internal rays
(inside basins of periodic attractors). This means in particular
 that any  critical point of $f$ which is attracted to a periodic attractor,
is outside the filled Julia set of the complex box mapping $F$. 
One of the reasons for setting-up the construction in this
way, is that this will allow us to show that 
puzzle pieces shrink to zero in diameter. 

\subsection{Puzzle construction for polynomials with 
periodic attractors or escaping critical points}
\label{subsec:puzzles}

The next construction coincides with the usual Yoccoz puzzle
construction in the case where the Julia set is connected
and there are no attracting fixed points. 
First choose a level curve $\Gamma$ of the
equipotential function  of the Julia set 
such that the following three properties hold for any
connected component $K'$ of the filled Julia set $K$:

1) if $K'$ contains a critical point and is periodic,
i.e. there exists some (minimal) $n>0$ so that $f^n(K')=K'$, then 
we want that each component of $\C\setminus \Gamma$
contains at most one of the sets
$K',\dots,f^{n-1}(K')$;

2) if $K'$ contains a critical point and
is not periodic, i.e. if $K',f(K'),\dots$ are all mutually
disjoint, then we want that not all these iterates of $K'$
are contained in one component of $\C\setminus \Gamma$;

3) if two critical points lie in one component of 
$\C\setminus \Gamma$
then they lie in the same component of $K$;

4) if a critical point
lies in a bounded component of $\C\setminus \Gamma$
then it lies in $K$.

Since equipotential level curves can be chosen to lie inside 
an arbitrary small neighbourhood of $K$ it is possible
to choose $\Gamma$ as above. 
Now take a component $K'$ of
$K$ which is periodic and contains a critical point;
assume its period is $n$ and  
let $d$ be the degree of $f^n|{\mathcal V}'$ where ${\mathcal V}'$
is the component of $\C\setminus \Gamma$ containing $K'$.
By \cite{LePrz}, each repelling
periodic point in $K'$ has at least one and at most a finite number
of external 
rays landing on it, and 
the number of rays of period $n$ landing on
$K'$ (i.e., on fixed points of $f^n|K'$)
is equal to the number of periodic points of
$\theta\mapsto d\, \theta  \,\, \mod 2\pi$. 
We are primarily interested in periodic points which 
are {\em separating}, i.e. such that there are several rays 
landing on it;  rays landing on such a periodic point
divide the complex plane in several components.

{\bf Case 1:}  Let us first assume that all fixed points of $f^n|K'$ are repelling. 
Then,  since $f^n|K'$ has $d$ (repelling) fixed points and 
$\alpha \mapsto d\alpha \,\, \mod 2\pi$ has only $d-1$ fixed points,
at least one of the fixed points is separating,
i.e., has more than one ray (with period $>n$) landing on it.

{\bf Case 2:}  If at least one of the periodic
 points of $f^n|K'$ is attracting, 
then the counting argument used in  case 1
might not work. By replacing $f^n$ by an iterate
we may assume that all periodic attracting points 
of $f^n\colon K'\to K'$ are fixed
points.
So consider the immediate basin $B$ of a hyperbolic
attracting  fixed point $p$ of $f^n|K'$. The following lemma
gives that there is a periodic curve in $B$ which  
has exactly period $2$ under $f^n$:

\begin{lemma} Let $g$ be a polynomial
with an attracting hyperbolic fixed point $p$ with
immediate basin $B$. Then there exist 
\begin{enumerate}
\item a curve
$\gamma$ of period $2$ in $B$ which connects a point of 
period precisely two in $\partial B$ to $p$ such that
$\gamma\cap g(\gamma)=\{p\}$ and such that
$\gamma\cup g(\gamma)$ divides $B$ into two components;
\item a curve $\Gamma_B\subset B$ surrounding $p$
such that $(g|B)^i(\Gamma_B)$ and 
$(g|B)^j(\Gamma_B)$ are disjoint for all $i,j\in \Z$, $i\ne j$.
\item $\gamma$ and $(g|B)^i(\Gamma_B)$ intersect in a single point,
for any $i\in \Z$.
\end{enumerate}
\end{lemma}
\begin{proof}
Taking the Riemann mapping $h\colon \D\to B$
we get that $h^{-1}\circ g\circ h$ is a proper mapping
of $\D$ onto itself and so a Blaschke product $\varphi$ with an attracting
fixed point $p_\varphi$ in $\D$. Note that all periodic orbits of
$\varphi|\partial \D$ are repelling. Hence there exists $N$ so that
$|D\varphi^N_{| \partial \D}|>2$. Let $V$ be a 
neighbourhood of $\partial \D$
such that $|D\varphi^N_{|V}|>2$ and such that $\varphi^{-N}(V)$ is compactly
contained in $V$.

Next choose a  
curve $\hat \gamma\subset \D$ ending at $p_\varphi$  such that 
$\varphi^2(\hat \gamma)\subset \hat \gamma$ and so that 
$\varphi(\hat \gamma)\cap \hat \gamma$ is equal to $p_\varphi$. 
This is easy to do:
First take a fundamental annulus $A$ around $p_\varphi$.
Then choose a curve $\tau$ in $A$ connecting some point $x$ in the outer boundary of $A$ with a point $y\ne \varphi(x)$ in the inner boundary of $A$ and next choose a curve  $\tau'$ in 
$\varphi(A)$ connecting $y$ to $\varphi^2(x)$ with the property that 
$\varphi(\tau)$ does not intersect $\tau'$.
The curve $\hat \gamma=\cup_{i\ge 0} \varphi^{2i}(\tau\cup \tau')$ has the required properties (usually $\hat \gamma$ will 
spiral towards $p_\varphi$). 
Now consider backward iterates of $\tau\cup \tau'$. 
Taking appropriate preimages (i.e. components) 
$\cup_{i} \varphi^{-2i}(\tau\cup \tau')$ becomes a curve
$\tilde \gamma$; we can assume that $\tau,\tau'$
are chosen so that
$\tilde \gamma$ is disjoint from the postcritical set. For $i$ sufficiently large $(\varphi|\D)^{-i}(\tau\cup \tau')$ lies
in $V$.  Hence, by the above expansion properties, the length of the curve $\tilde \gamma$ is finite.
Since $\varphi^2(\tilde \gamma)=\tilde \gamma$, it follows that 
$\tilde \gamma$ approaches a point $q\in \partial \D$ which has precisely period two
under $\varphi$. Moreover, $\tilde \gamma$ approaches $q$
non-tangentially: this holds because $\varphi$ is expanding
near $q$ and leaves invariant $\partial \D$. Hence, by  
\cite[Theorem 1]{Pom}, the corresponding curve 
$\gamma:=h(\tilde \gamma)$ lands on a fixed point of 
$g^2$. In fact, more holds.
{\em Claim:}  $\gamma$ lands on a periodic point
of $g$ of period exactly  two. 
Indeed, assume by contradiction that  $\gamma$ lands on a fixed point. Then let $B_+$ one of the components
of $D\setminus (\gamma\cup g(\gamma))$.
Since the boundary of $B_+$ is equal to the Jordan curve
$\gamma\cup g(\gamma)$
and since $g^2(\gamma)=\gamma$,
the boundary of the set $g(B_+)$ again  consists of the Jordan curve
$\gamma\cup g(\gamma)$.
Since $g(B_+)\ne B_+$ it follows from the Jordan curve
theorem that $g(B_+)$ is the unbounded complement of 
$\gamma\cup g(\gamma)$. This contradicts the assumption that
$g(B_+)\subset D$, and finished the proof of the claim.

To complete the proof of the lemma, 
we notice that the curve $\Gamma_B$ can be constructed using a preimage of 
the outer boundary of $A$ (which can be constructed by considering
the boundary of a neighbourhood of $p$ on which $g$ is linearizable).
We can and will choose this curve so that it does not contain iterates of critical points.
\end{proof}

The curve $\gamma$ ('internal ray') 
together with the external rays landing 
on the corresponding period two points 
of $f^n|\partial B$, are called the {\em rays
associated to the attracting fixed point} of $f^n|K'$.
These curves are smooth except at the attracting fixed point
and its (other) endpoints.
Of course the union of such internal and external 
rays again separates the plane.
The curve $\Gamma_B$ is called an {\em equipotential} 
of the attracting
fixed point. We can choose $\Gamma_B$ so that 
it bounds a disc containing $p$ and all iterates of critical points of $f$
which lie in $B$.

For each component $K'$ of the filled Julia
set containing a critical point, denote its period by $n=n(K')$, 
and consider the union of rays associated to the 
attracting fixed points of $f^n|K'$ together with all
external rays  landing on separating fixed points of
$f^n|K'$. 
Moreover, consider the equipotentials $\Gamma$ and
the union of all 
equipotentials $\Gamma_B$ associated to attracting orbits; 
the set ${\mathcal V}$ is defined to be the union of  the bounded components of $\C\setminus \Gamma$ and minus
the bounded components of $\C\setminus \cup \Gamma'$
where the last union runs over all $\Gamma'$ as in Case 2.
Let the partition $\Part=\{P_1,\dots,P_k\}$ of 
${\mathcal V}$ be defined by all such rays, i.e. as the components of $\mathcal V$ minus the above collection of rays.
Combining Cases 1 and 2 above, we have shown that
if a  component of ${\mathcal V}$ contains a critical point then it
is partitioned in two or more pieces from $\Part$.
(If there are no components of $K$ 
containing critical points
then the partition does not involve rays.)

Components of $f^{-j}(\Part)$ map onto a union of 
components of $f^{-(j-1)}(\Part)$.
Components of $f^{-j}(\Part)$ are called {\em puzzle pieces}
and the pieces which contain
critical points are called {\em critical puzzle pieces}.

Note that  each puzzle piece is the closure of a finite union of open arcs,
consisting of equipotentials and rays.

\subsection{How to associate a complex box mapping to a  polynomial}
\label{subsec:boxmaps}
For critical points $c,c'$ we define 
$c \preceq c'$ if the forward orbit of $f(c')$ intersects
each puzzle piece containing $c$. 
If $c\preceq c$ then $c$ is called {\em puzzle recurrent}.
The partial ordering $\le $ is defined by
$c\le c'$  if $c=c'$ or  $c\preceq c'$.

\begin{lemma}\label{lem:constrbox}
Let $f$ be a  finitely renormalizable polynomial $f$
with all periodic points hyperbolic, and consider the puzzle
construction from above. Then for
any puzzle recurrent critical point $c$,
there is a critical puzzle piece $W$ 
such that
\begin{enumerate}
\item the component of the first return map to $W$ containing $c$
is compactly contained in $W$;
\item the first return map to $W$ is non-renormalizable;
\item take $V$ to be the union of the 
components of the first entry map to $W$
which contain a critical point, and $U\subset V$ the 
domain of the first return  map to $V$.
Then $R\colon U\to V$
is a complex box mapping
with the property that $V$ contains 
no escaping critical point of $f$.
\end{enumerate}
\end{lemma}
\begin{proof} The proof we give here simplifies an earlier
proof, and is  due to Genadi Levin. 
First assume that $f$ is non-renormalizable.
We claim that for each critical point $c$
there is a critical piece $P\ni c$, so that its
closure is disjoint from all separating periodic points used
for the construction of the puzzle partition ${\mathcal P}$.
Indeed, otherwise all critical pieces containing $c$ have
a (separating) periodic point $a$ in their boundary, and
hence the intersection of (all) the critical pieces of $c$
is a continuum which also contains $a$.
Because $f$ is non-renormalizable 
this continuum must be wandering, see \cite{BlokhLevin}.
But since the continuum contains the fixed point
$a$, this is impossible.

Let us fix $P_0=P\ni c$. Since the closure of
$P$ contains no fixed point, and the boundary of $P$
consists of pieces of equipotentials and rays, 
the boundary of $P$ contains no periodic points.
Now define inductively $P_{i+1}$ to 
be the first return domain to $P_i$ containing
$c$.  Assume that $P_{n+1}$ is not compactly contained
in $P_n$. Then there are two points
$p_1, p_2$ in the boundary of $P_n$, such that
$f^k(p_1)=p_2$, where $k$ is the return time
of $P_{n+1}$ to $P_n$. 
Let $r$ be so that $f^r(P_n)=P_0$, and 
$q_1=f^r(p_1), q_2=f^r(p_2)$.
Then $q_1, q_2$ are in the boundary of $P_0$,
and $f^k(q_1)=f^r(f^k(p_1))=f^r(p_2)=q_2$.
If $n$ is large, then $k$ is large too ($f$ is non-renormalizable), 
and $f^k(q_1)=q_2$
is impossible, because the boundary of $P_0$ consists of
rays and equipotentials (not depending on $k$) 
and its boundary is disjoint from periodic points.
Taking $W=P_n$ with $n$ minimal as above
the lemma follows. 

Now assume $f$ is renormalizable. In this case 
we have to work with  rays associated
to generalized polynomial-like maps (here we allow
that the domain $U$ of polynomial-like maps $F\colon U\to V$
we consider has infinitely many components) .
Let $P_0$ be as in the non-renormalizable case 
and define  inductively  $P_{i+1}$ to be the 
component containing $c$ of the  domain of definition of
the  first return domain map $R_i\colon \dom(P_i)\to P_i$.
If for all $i$, there exists some
$k>0$ so that $R_i^k(c)\notin P_{i+1}$ 
then $R_i\colon \dom(P_i)\to P_i$
is non-renormalizable at $c$ and we are in Case 1.
So we may assume that there exists a level $i$ so that $R_i^k(c)\in P_i$ for all $k\ge 0$.
In fact, we can take $i$ so large that
each other critical point $c'\ne c$ in $P_{i}$ also does not leave
$P_i$ under iterates of $R_i\colon \dom(P_i)\to P_i$.
Then the filled Julia set $K_i$ of $R_i\colon \dom(P_i)\to P_i$ 
is connected and not a singleton.
From Theorem $1^\tau$ of \cite{LePrz} each fixed point 
of $R_i\colon \dom(P_i)\to P_i$  
has external rays landing on it, and there is at least one
such fixed point with a ray of some period $q>1$ landing on it.
Let $\gamma$ be the union of rays of fixed points of 
$R_i\colon \dom(P_i)\to P_i$
with rays of period $>1$ landing on it, and consider 
$P_i\setminus \gamma$. Since the rays have period $>1$, 
$c$ and $R_i(c)$ cannot lie in the same component of 
$P_i\setminus \gamma$. Let $P^1$ be the component
of $P_i\setminus \gamma$ containing $c$, and consider
the first return map $R\colon \dom(P^1)\to P^1$. 
Now repeat the above  construction, replacing $R\colon \dom(P)\to 
P$
by $R\colon \dom(P^1)\to P^1$.  Since $f$ is finitely renormalizable
(say of period $s$) this cannot be repeated more than $k$ times
where $2^k\le s$. 
\end{proof}

\noindent
{\bf Question:} take a polynomial with all periodic points hyperbolic; 
is it always possible to find puzzle pieces which contain at most one critical point?  

We are only able to answer this question once we have  complex bounds (so only in the 
case when the polynomial is only finitely renormalizable). That is the reason why in the 
definition of complex box mappings we allow several critical points to be in one component
of the domain.

\bigskip
Now we will show how to use the previous lemma
to construct a complex box mapping.

\begin{corollary}
  \label{cor:boxmap} For a polynomial as in the previous lemma, there
  exists a collection of puzzle pieces $V$ which contains all critical
  points which do not escape to infinity or are in the basin of
  periodic attractors (and only such critical points), such that the
  first return map into $V$ defines a complex box mapping $F\colon U
  \to V$.
\end{corollary}
\begin{proof} Take a critical point $c_1$ which is puzzle recurrent
and consider the critical puzzle piece $W_1\ni c_1$  from the previous lemma (if there is no such critical point
then skip this step), and write ${\mathcal W}_1=W_1$.  
Now continue by induction and assume ${\mathcal W}_{k-1}$ is chosen.
Then choose a critical point $c_k$ which is 
puzzle recurrent and which is never mapped into ${\mathcal W}_{k-1}$. (If this is impossible, the induction stops, we set
$m=k-1$ and this part of the construction is completed.)
Then let $W_k\ni c_k$ be the critical puzzle piece from the previous lemma and whose level is at least the level of each of 
the puzzle pieces in  $\mathcal W_{k-1}$. Then 
write ${\mathcal W}_k={\mathcal W}_{k-1}\cup W_k$. After a finite number of steps,  the  only critical points which are 
not contained in or eventually mapped into some union of puzzle pieces ${\mathcal W}_m$, are those which are not puzzle recurrent.
(If there are no puzzle recurrent critical points, then define 
$m=0$ and ${\mathcal W}_0=\emptyset$).

Next choose a critical point $c_{m+1}$ which is not puzzle recurrent
and which is never mapped into ${\mathcal W}_m$, choose a puzzle piece $W_{m+1}\ni c_{m+1}$ of level at least the level of each of the puzzle pieces of ${\mathcal W}_m$ and such that $c_{m+1}$ never
re-enters 
$W_{m+1}$. Then write ${\mathcal W}_{m+1}={\mathcal W}_m\cup W_{m+1}$. Continue in this way as long as possible (for $m'$ steps).
Next let $V$ be the union of ${\mathcal W}_{m'}$ and
the domains of the first return map to ${\mathcal W}_{m'}$
containing critical points outside ${\mathcal W}_{m'}$. 
The set $V$ consists of puzzle pieces  and
contains all non-escaping 
critical points.  $V$ is {\em nice}:
the boundary of $V$ is never mapped into the closure
of $V$.
The return map to $V$ therefore defines a box mapping, whose
domain $U$ is equal to the set of points in $V$ 
which eventually return  to $V$.
\end{proof}

A critical point which is Misiurewicz (say eventually mapped
to a repelling periodic orbit) will be contained 
in a component of   $V$, but might not be in $U$.

\bigskip

Now we will finally show that Theorem~\ref{thm:main2}
implies Theorem~\ref{thm:main}.
Indeed, Theorem~\ref{thm:main2} implies that all puzzle pieces
of this complex box mapping shrink to zero, and that it carries no
measurable invariant  line field on its Julia set. The set of points
which are never mapped into $V$ (and is not contained
in the basin of a periodic attractor) is hyperbolic, and therefore
puzzle pieces containing such points also shrink to zero.
Theorem~\ref{thm:main} therefore follows.
\bigskip

\section{QC-conjugacies on basins of periodic attractors}

\begin{lemma}\label{lem:basins}
Assume that $f$ and $\tilde f$
are finitely renormalizable and only have hyperbolic periodic points. 
\begin{enumerate}
\item If $f$ and $\tilde f$ are topologically conjugate on a periodic 
component of the Fatou set, then this conjugacy extends
continuously to the boundary of this component;
\item If $h$ is a topological conjugacy between $f$ and $\tilde f$
 restricted to their Fatou sets, then $h$ can be replaced by 
a qc conjugacy $h_0$ such that $h$ and $h_0$  have the same extension to the boundary of periodic components of the Fatou set.
(In particular, the extensions agree on the boundary of the basin
of $\infty$, i.e. on the Julia set.)
\end{enumerate}
\end{lemma}

\noindent
Remark that inside the level set associated to the B\"ottcher
coordinates of a super-attracting fixed point there exists
a dense set $X$ so that for each $a,b\in X$ there exists $n>0$
so that $f^n(a)=f^n(b)$. 
(So the grand orbit of a point in the basin of a super-attracting fixed point
is dense in such a topological circle.)   It follows that a topological conjugacy
maps level sets of B\"ottcher coordinates to level sets.
The same
argument also shows that the conjugacy maps
rays to rays.

\begin{proof} 
First consider a periodic component $B\subset \bar \C$ 
of the Fatou set
which contains a super-attracting periodic point (of, say, period $n$).
Consider an equipotential $\Gamma'$ 
on which $f^n$ is conjugate to $z\mapsto z^d$, $z\in \D$.
(Choose $\Gamma'$ so that no iterates of
critical points of $f$ are on $\Gamma'$, and so that
the disc which $\Gamma'$ bounds contains no other
critical point apart from the fixed point.)
The curves $\Gamma'$ and 
$f^n(\Gamma')$  bound a fundamental annulus (which each orbit
hits at most twice).
By the above remark, the conjugacy $h$ maps this annulus to a similar annulus again bounded by equipotentials, 
and  so $h$ can be written in polar B\"ottcher coordinates 
as 
$$(r,\phi) \mapsto (Z(r),\phi).$$
Now approximate $Z$ by a smooth map $Z_0$
which agrees with $Z$ 
on the inner and outer boundary of the annulus
and also on circles containing iterates of critical points.
By its form $(r,\phi)\mapsto (Z_0(r),\phi)$ induces 
a qc conjugacy $h_0$ between $f$ and $\tilde f$ on the 
the disc in $B$ bounded by $\Gamma'$. 
Since $h_0$ agrees with $h$
on iterates of critical points of $f$, one can extend
$h_0$ to $B$ by pulling back. Similarly, 
near an attracting periodic
point of period $n$ and with multiplier $\ne 0$, one can do something similar. Indeed, in this case
take a curve $\Gamma'$ surrounding the periodic point,
which in linearizing coordinates is a circle, 
and such that there are no iterates of critical points on $\Gamma'$.
Also choose $\Gamma'$ so that $f^n$ is univalent on the disc
bounded by $\Gamma'$. Then approximate the conjugacy $h$ on the (fundamental) annulus
bounded by $\Gamma'$ and $f^n(\Gamma')$ by a smooth
homeomorphism $h_0$ which agrees with $h$ at iterates of
critical points of $f$. Do this, so that $\tilde \Gamma'=h_0(\Gamma')$
is smooth, and so that $h_0$ maps $f^n(\Gamma')$ to 
$\tilde f^n(\tilde \Gamma')$.
Again this induces a conjugacy $h_0$ between
$f$ and $\tilde f$ on the basin of this periodic attractor.

Let us show that $h_0$ (or $h$) can be extended continuously 
to the boundary of $B$. So consider the immediate
basin $B$ of an  attracting periodic point of period $n$.
Let $\phi\colon \D\to B$ be the Riemann mapping. Then
$A=\phi^{-1}\circ f^n \circ \phi$ is a Blaschke product
(say of degree $d$)
and since all periodic points of $f$ are repelling,  
$A$ is expanding on a neighbourhood $V$ of $\partial \D$
(i.e. there exists $N$ such that  $A^{N}(V)\supset V$
and $|DA^N|\ge 2$ on $V$). By taking components $\Gamma_i'$
of $f^{-ni}(\Gamma')$ such that $\Gamma_0',\Gamma_1',\dots$
lie nested, we obtain a nested sequence of fundamental annuli 
$Fund_i$  bounded by $\Gamma_i'$ and $\Gamma_{i+1}'$.
For $i$ large enough,
$\phi^{-1}(Fund_i)$  is contained in $V\cap \D$. For simplicity write
$Fund:=Fund_i$.  It follows that if we consider a smooth curve $\tau$
connecting a point $z$ in the outer boundary of $Fund$ to 
$f^n(z)$, then the length of each component 
of $A^{-i} (\phi^{-1}(\tau))$ decreases exponentially
with $i$. Hence $\cup_{i\ge 0} A^{-i}(\phi^{-1}(\tau))$ contains
a curve $\tau'$  which converges to a fixed point of 
$A|\partial \D$. Because of this, 
$A^{-i}(\phi^{-1}(Fund\setminus \tau'))$ consists of
$d^i$ topological rectangles with two of its boundaries consisting of
curves through preimages of one of the fixed points of $A|\partial \D$
(these two curves are components of   $A^{-i}(\tau')$).
In other words, each of the $d^i$ preimages of $A^{-i}(z)$
with $z\in \phi^{-1}(Fund)$ is determined by a sequence of $i$
symbols. Now do the same for $\tilde f$. Because of the symbolic
description given above, and since the size of each of the $d^i$
rectangles shrinks to zero as $i\to \infty$, it follows that
$\phi^{-1}  \circ h_0 \circ \phi $ extends continuously to a 
conjugacy on the boundary. 
But since $J(f)$ and $J(\tilde f)$ are locally 
connected, it follows that $h_0$ extends continuously
to the closure of $B$. 

Of course such an extension is not unique
if the degree of $A$ is $>2$.
But note that $h_0|\partial B$ agrees 
with $h|\partial B$ because $h$ and $h_0$
agree on each of the points
$\phi(\tau')\cap \partial f^{-i}(Fund)$, $i=0,1,\dots$, 
i.e. on a sequence
of points in $\phi(\tau')$ converging to $\partial B$.

By doing this for each periodic attractor, and pulling back
the resulting map $h_0$ we obtain the required homeomorphism.
\end{proof}

\section{Rigidity of non-renormalizalble polynomials (Theorem~\ref{thm:main2} implies Theorem~\ref{thm:rigidity})}

Let us now show that combinatorially equivalent
polynomials induce combinatorially equivalent 
complex box mappings.
\medskip

\begin{proposition}\label{prop:combequiv}
Assume that $f$ and $\tilde f$ are finitely renormalizable polynomials with only hyperbolic periodic points. 
Moreover, assume that one of the following
conditions hold:
\begin{description}
\item[a)\quad] $f$ and $\tilde f$ are
topologically conjugate;
\item[b)\quad] $f,\tilde f$ are combinatorially equivalent 
and restricted  to their Fatou sets there exists a
compatible conjugacy between  $f$ and $\tilde f$
(these notions are defined in Definition~\ref{defcombequiv}). 
\end{description}
Then associated to $f,\tilde f$ are complex
box mappings $F\colon U\to V$, $\tilde F\colon \tilde U \to \tilde V$
which are combinatorially equivalent
and there exists a qc map $h$ which maps $V$ onto $\tilde V$
and for which $\tilde F\circ h = h \circ F$ on $\partial U$.

\end{proposition}

In other words, 
the conclusion is that the 
assumption of Theorem~\ref{thm:main2} is satisfied.
As a first step in the proof of this proposition we have:

\begin{lemma} Assume $f,\tilde f$ are as in the proposition above.
Then, associated to $f$ and $\tilde f$
are complex box mappings $F\colon U\to V$ and $\tilde F\colon \tilde U\to \tilde V$ so that to each puzzle piece $P$ of $F$ is uniquely associated to a puzzle piece $\tilde P$ of $\tilde F$ with the following properties:
\begin{enumerate}
\item $P$ contains the same number of
 critical points as $\tilde P$;
\item the conjugacy from Assumption 1 or 2 of the above proposition, 
maps $P\setminus J(f)$ onto $\tilde P\setminus J(\tilde f)$;
\item $\widetilde{f^k(P)}=\tilde f^k(\tilde P)$.
\end{enumerate}
\end{lemma}
\begin{proof} 
Choose equipotentials $\Gamma$ and $\Gamma_{BS}$ (associated
to $\infty$ and super-attracting periodic orbits in the previous subsection) which do not  contain
iterates of critical points of $f$. Then $h(\Gamma)$
and $h(\Gamma_{BS})$ are also equipotentials, see 
the remark above, and so define $\tilde \Gamma=
h(\Gamma)$ and $\tilde \Gamma_{BS}=h(\tilde \Gamma_{BS})$.
For periodic attractors, choose an equipotential $\Gamma_B$ as in 
Subsection~\ref{subsec:puzzles} and set $\tilde \Gamma_B=h(\Gamma_B)$. Similarly, choose rays for $f$ which are mapped by 
$h$ onto rays. Using these choices,
associate to $f$ and $\tilde f$
complex box mappings $F$ and $\tilde F$
as in Subsection~\ref{subsec:boxmaps}
(using  the same construction for $f$
as for $\tilde f$).  

For $\alpha\in \Q$, let $p(\alpha,f)$ be the periodic point on which the $f$-ray with angle $\alpha$ lands and define $n(\alpha,f)$ to be the number of rays landing on $p(\alpha,f)$.
By definition,  $n(\alpha,f)=n(\alpha,\tilde f)$. 
So if we determine a puzzle piece by the rays 
in its boundary, then a collection of rays
determines a puzzle piece for $F$ if and only if it determines a puzzle
piece for $\tilde F$. A puzzle piece contains a critical point if and only if
two of its boundary rays are mapped onto a single ray.
Hence the lemma follows. 
\end{proof}

\begin{proofof}{Proposition~\ref{prop:combequiv}}
Let us now show that $F$ and $\tilde F$ are combinatorially 
equivalent (in the sense of Definition~\ref{def:equivbox})
and that the assumption of Theorem~\ref{thm:main2} holds.
Let $U_0,\tilde U_0$ be  collections
of corresponding puzzle pieces of $F,\tilde F$ 
containing the set of critical points of $F,\tilde F$
(i.e., the component of $U,\tilde U$ containing critical points
of $F,\tilde F$).

By Lemma~\ref{lem:basins}, we can  assume that
the (compatible) conjugacy on the Fatou set, induces 
a $K_0$-qc homeomorphism $h\colon V\to \tilde V$
mapping $U$ onto $\tilde U$
such that $\tilde F \circ h = h \circ F$
and which is a combinatorial equivalence (in the sense
of Definition~\ref{def:equivbox}). Of course, $h$ induces
a $K_0$-qc map
 $h\colon V\setminus K(F) \to \tilde V \setminus K(\tilde F)$
 which conjugates $F$ and $\tilde F$ on these sets.

Let us now cite the Spreading principle (see Section 5.3 of \cite{KSS1}):

\newtheorem*{spreading}{Spreading Principle}
\begin{spreading}
In the above setting, 
there exists a $K_0$-qc map $\Phi:\mathbb{C}\to
\mathbb{C}$
such that the following hold:
\begin{enumerate}
\item $\Phi=\phi$ on $U_0$;
\item for each $z\not\in U_0$ we have
$\tilde{f}\circ \Phi(z)=\Phi\circ f(z)$;
\item The dilatation of $\Psi$ on $\mathbb{C}-D(U_0)$ is equal to $K_0$, where $D(U_0)$
denotes the domain of the first landing map under $f$ to $U_0$;
\item  for each puzzle piece $P$ of $F\colon U\to V$
which is not contained in $D(U_0)$, $\Phi(P)$ is equal to the 
puzzle piece
$\tilde{P}$ of $\tilde F\colon \tilde U\to \tilde V$ and
$\Phi:P\to {\tilde{P}}$
agrees with $h$ their boundary.
\end{enumerate}
\end{spreading}

\begin{proof} 
Denote by $\Part_n$ the collection of level $n$ puzzle pieces.
For a puzzle piece $P\in \Part_n$, let $k=k(P)\leq n$ be the minimal
non-negative integer such that $F^k(P)$ is a critical puzzle piece or
has depth $0$, and let $\tau(P)=F^k(P)$. Then $F^k: P\to \tau(P)$ is a
conformal map, and so is $\tilde{F}^k:\tilde{P}\to
{\tau(\tilde P)}$.
Now define a qc map
$\phi_P:P\to \tilde{P}$ by the formula $\tilde{F}^k\circ \phi_P=
h\circ F^k$,
where $h$ is the $K_0$-qc map from above. Note that
the map $\phi_P$ has the same maximal dilatation as $h$.

Let $Y_0$ be the union of all puzzle pieces in $\Part_0$.
For $n\geq 0$, inductively define $Y_{n+1}$
to be the subset of $Y_n$
consisting of puzzle pieces $P$ of depth $n+1$ so that
$P$ is not contained in $D(U_0)$, the domain of the first entry map
into $U_0$.
Note that each puzzle piece in $Y_n-Y_{n+1}$ of depth
$n+1$ is a component of $D(U_0)$.

Define $\Phi_0=h$.  For each $n\geq 0$, assume that
$\Phi_n$ is defined and define $\Phi_{n+1}$ so that
$\Phi_{n+1}=\Phi_n$ on $V-Y_n$ and 
for each
component $P$ of $Y_n$ define  
\begin{enumerate}
\item $\Phi_{n+1}=h$ on $P-\bigcup_{Q\in\Part_{n+1}}Q$;
\item  for each $Q\in\Part_{n+1}$ which is
contained in $P$, define 
$\Phi_{n+1}=\phi_Q$.
\end{enumerate}
For each $n\geq 0$, $\Phi_{n}$ is a $K_0$-qc map. Note that $\Phi_n$
is eventually constant (not dependent on $n$ for large $n$) 
on $\mathbb{C}-\bigcap_n Y_n$. Since $\bigcap Y_n=E(U_0)$ is
a nowhere dense set, $\Phi_n$ converges
to a qc map $\Phi$. The properties (1), (2) and (4) follow directly
from the construction, and (3) follows from the fact that $E(U_0)$
(which consists of points which are 
never mapped into $U_0$) has measure
zero.
\end{proof}

Let us now complete the proof of Proposition~\ref{prop:combequiv}.
Taking for $U_0,\tilde U_0$ the collection of critical puzzle piece of level  $n$, set $H_n=\Psi$.
Then $H_n\circ F=\tilde F\circ H_n$ on $V\setminus F^{-n}(V)$.
Moreover, $H_n$ agrees with the compatible conjugacy $H$
on $V\setminus F^{-1}(V)$.
Since $H$ is $\max(K,K_0)$-qc, the proposition follows.
\end{proofof}

It follows that  Theorem~\ref{thm:main2} implies 
Theorem~\ref{thm:rigidity}.

\section{Robustness of complex box mappings}
\label{sec:robust}

Let $V$ be the set consisting
of a union of puzzle pieces which contains each critical point which
is not in the basin of a periodic attractor and 
let
$$E(V,f)=\{z\in J(f) ; f^n(z)\notin V\mbox{ for all }n\ge 0\}.$$
For each $z\in E(V,f)$ the level $k$ puzzle pieces $P_k(z)$ containing $z$ 
forms a decreasing sequence shrinking  to $z$, see Lemma
\cite[Lemma 1.8]{MR1765085}. It follows that $E(V,f)$ is hyperbolic.

 So if $f$ has only hyperbolic periodic points, there exists 
 a neighbourhood $O$ of $f$ so that 
 if $V_g$ is the puzzle piece corresponding $V$ 
 for $g\in O$  then $E(V_g,g)$ is hyperbolic
 and for each $z\in E(V,f)$
 these exists an analytic map 
 $$O\ni g \mapsto z(g) \in E(V_g,g).$$

Take a polynomial $f$ with all periodic points hyperbolic, 
with critical points $c_1,\dots,c_d$ and with
complex box mapping $F\colon U\to V$ 
as in Subsection~\ref{subsec:boxmaps}.
Then there exists (locally) a manifold $\hat \Sigma(f)$ 
containing $f$  (also depending on $V$)
with the following properties: each $g\in \hat \Sigma(f)$ has critical points 
$c_1(g),\dots, c_d(g)$
such that
\begin{enumerate}
\item 
$g\mapsto c_i(g) $ is continuous
and  the degree of $g$ at  $c_i(g)$ 
is the same for each $g\in \hat \Sigma(f)$;
\item let $\{A_i\}$ resp.  $\{SA_i\}$ be the collection of hyperbolic attracting (resp. super-attracting) periodic points of $f$;
then choose $O$ small enough so that 
$O\ni g\mapsto A_i(g)$ and  $O\ni g\mapsto SA_i(g)$ 
are continuous (but $g$ might have new periodic attracting orbits);
\item if $f^N(c_i(f))=f^M(c_j(f))$ for some  integers $i,j$ and $N,M$ and this point is
in the basin of a hyperbolic periodic attractor,  then for all $g\in \hat \Sigma(f)$ 
we have $g^N(c_i(g))=g^M(c_j(g))$; by decreasing $O$ if necessary we
can ensure that if 
for some  integers $N,M,i,j,k$ and for some $g \in \hat \Sigma(f)$, \,\, 
 $g^N(c_i(f))=g^M(c_j(f))$ and this point is in the
basin of a hyperbolic attractor $A_k(g)$ (which is a continuous
deformation of an attractor of $f$) , then  $f^N(c_i(f))=f^M(c_j(f))$;
\item if  $c_i(f)$ is in the basin of a super-attracting fixed point
with the angle of the ray landing on $c_i(f)$ equal to $\alpha_i$
then the same holds for all 
$g\in \hat \Sigma(f)$ (so the angle is the same);
\item if for some integers $i,j$ and $N,M$, \,\,
$f^N(c_i(f)),f^M(c_j(f))$ are both in the basin of a super-attracting fixed point while hitting 
the same  level set associated to its B\"ottcher coordinates then
the corresponding statement  holds for all $g\in \hat \Sigma(f)$;
\item if the $k$-th iterate  of $c_i(f)$ is equal to $z\in E(V)$
then $g^k(c_i(g))=z(g)$ where $z(g)$ is as above.
\end{enumerate}
Using quasi-conformal surgery and transversality of
unfoldings  (see for example \cite{StrFam}) it is not hard to show that 
$\hat \Sigma(f)$ is a manifold (locally near $f$).

\begin{proposition}\label{prop:nearbypol}
Let $f$ be a non-renormalizable 
polynomial with all periodic points hyperbolic, 
with critical points $c_1,\dots,c_d$ and
let  $F\colon U\to V$ be a complex box mapping
associated to $f$ as in Subsection~\ref{subsec:puzzles}.
Then there exists a neighbourhood $O$ of $f$
with the following properties. If $\tilde f\in O\cap \hat \Sigma(f)$ then 
there exists a complex box mapping $\tilde F\colon \tilde U\to \tilde V$ associated to $\tilde f$
and a qc map $h\colon \C\to \C$ which maps
$U$ to $\tilde U$ and such that $\tilde F\circ h =h\circ F$ on 
$\partial U$.
\end{proposition}

\begin{proofof}{Proposition~\ref{prop:nearbypol}}
 Let $P_0$ be the puzzle piece partition for $f$ of level $0$
  constructed using the external rays landing on the hyperbolic
  periodic point $a$  and an equipotential level.
 For all
  polynomials $g$ close to $f$, the number of rays landing 
  on $a(g)$ is the same. So one can define a puzzle piece $P_0(g)$
  corresponding to $P_0$ in an obvious way (taking the rays landing
  on $a(f)$ and the equipotential of  the same level). 
  Clearly, each puzzle piece in $P_0(g)$ converges to a
  puzzle piece of $P_0$ when $g\to f$ (say in the Hausdorff topology). 
Now take $\tilde f$ so close to $f$ that
each periodic point used in the construction of the complex
box mapping $F\colon U\to V$ remains hyperbolic.
It then follows that there exists a finite collection
of puzzle pieces $\tilde V$ corresponding to $V$, and that the first return map
to $\tilde V$ induces a complex box mapping 
$\tilde F\colon \tilde U\to \tilde V$ corresponding
to $F\colon U\to V$.
The set of points which never enter $U$
is hyperbolic. Hence
since $\tilde f$ is close to $f$, any non-hyperbolic 
or 'new' attracting periodic points of $\tilde f$ is 
inside $\tilde U$.

Because of Assumptions 2 and 3 in the space $\hat \Sigma(f)$, 
we can construct a conjugacy between $f$ and $\tilde f$
on a neighbourhood of super-attracting periodic points
of the form $(r,\phi)\mapsto (Z(r),\phi)$ which maps iterates
of critical points to corresponding iterates of 
critical points. 
Similarly, one can construct a conjugacy near other periodic attractors
which are outside $V$. Using this, and the hyperbolicity of the
set $E(V)$, we obtain 
a qc conjugacy on the complement of $V$.
Using the spreading principle from the previous section
the result follows.
\end{proofof}

\section{Perturbations to hyperbolic
polynomials \\(Theorem~\ref{thm:main2} implies Theorem~\ref{thm:density})}

Let us say that a polynomial is {\em semi-hyperbolic} if all periodic points
are hyperbolic, and all critical points
in the Julia set are eventually mapped into some hyperbolic set.
By transversality of the unfolding of such a semi-hyperbolic map (see for example
\cite{StrFam}) it follows that it  can be approximated by   hyperbolic maps.

Let us consider the space ${\mathcal P}_d=\{z^d+a_2z^{d-2}+\dots+a_d\}$ of normalized polynomials  polynomials of degree $d>1$; 
if $f,\tilde f\in {\mathcal P}_d$ are conjugate
by a M{\"o}bius transformation they are equal, and every polynomial
is M{\"o}bius conjugate to one in ${\mathcal P}_d$.

  Let $f\in {\mathcal P}_d$ be a non-renormalizable polynomial of degree $d$ without neutral periodic points, and  assume for the moment  that all critical points of $f$ are non-degenerate. 
The map $f$ can have some critical points whose iterates converge     to periodic attractors. Since by Theorem~\ref{thm:main}
$f$ has no measurable invariant linefields on its Julia set,
by Theorem 6.3 in \cite{MS3} the space 
$QC(f)=\{g\in {\mathcal P}_d; g \mbox{ is qc conjugate to }f\}$ 
is a countable union of embedded manifolds $T_i$ each of 
at most complex dimension $s$, where $s$ is the
  number of (finite) critical points of $f$ which are in basins of periodic attractors. 
  Of course, $QC(f)$ is contained in the space
$\hat \Sigma(f)$ introduced in the previous section.
To express the dimension of 
$QC(f)$, we need to look at critical relations:
it is possible that there are 
critical points $c_i,c_j$ of periodic attractors which are 
critically related, because either $c_i=c_j$ or because 
\begin{enumerate}
\item some iterates of
$c_i$ and $c_j$ in the basin of a periodic attractor could 
be the same. In this case we say that $c_i$ and
$c_j$ are critically related; any map conjugate would need
to have the same properties;
\item some iterates of $c_i,c_j$
in the basin of a super-attracting periodic attractor
could hit the same level set associated to
the B\"ottcher coordinate. In this case we again say 
$c_i$ and $c_j$ are critically related. Any
conjugate  map would also have this property and the `external angles'
of these iterates of $c_i$ and $c_j$ would have to be the same.
\end{enumerate}
By the description in
Theorem 6.3 in \cite{MS3} the dimension of $QC(f)$
is equal to $s'$, where $s'$ is the number of equivalence classes 
of critical points in basins of hyperbolic periodic attractors.
Moreover, $\dim(\hat \Sigma)$ is the number of equivalence classes of
all critical points minus the number of critical points in the Julia set
which escape $V$ (see point 4 in the definition of $\hat \Sigma$).
So the codimension of $QC(f)$ in $\hat \Sigma(f)$ is equal to
the number of critical points in the Julia set which do not leave $V$.
 
Near $f$, the set $\hat \Sigma(f)$ is a smooth manifold 
and so take local 
coordinates in which $\hat\Sigma:=\hat \Sigma(f)$ is a linear space.
If $s'=\dim(\hat \Sigma)$ then all the critical points
which are not contained in the basin of a periodic attractor
are mapped into a hyperbolic set, and so the map is 
semi-hyperbolic and can be perturbed to a hyperbolic map
as observed in the beginning of this section. 
So assume $s'<\dim(\hat \Sigma)$ and 
consider the Grassman space $G$ of 
all complex linear subspaces $\Sigma\subset \hat \Sigma$ 
of complex codimension $s'$ 
containing $f$
(here $\Sigma$ is assumed to be a linear subspace in terms of the 
local coordinates which make $\hat \Sigma$  a linear space).
Denote by $B(f;r)$ the {\em open} ball  of radius $r$ through $f$. 
Since $\Sigma,T_i$ are both subsets of $\hat \Sigma$ 
such that the sum of their 
dimensions is $\dim(\hat \Sigma)$, for each $i$, the set of 
$\Sigma\in G$ and $r>0$ such that $\Sigma
\cap \partial B(f;r)\cap T_i=\emptyset$ 
is open and dense in $G\times \R^+$. By the Baire category theorem it follows that
there exist (many) $\Sigma_0\in G$ and $r_0>0$ such that 
$\Sigma_0\cap \partial B(f;r_0) \cap QC(f)=\emptyset$.
So throughout the remainder restrict ${\mathcal P}_d$ to $\Sigma_0$ and define 
$$W:=\Sigma_0\cap B(f;2r_0) \mbox{ and }W_0:=\Sigma_0\cap \overline{B(f;r_0)}$$
where we take $r_0$ so small that $W\subset \Sigma_0$
(we can still shrink $r_0$ later in the proof).
Let us denote the boundary of $W_0$ as a subset of $\Sigma_0$
by $\partial W$.
One reason for choosing  $W$ in this way is that  
$\partial W_0\cap QC(f)=\emptyset$.

We may assume that  $W$ is chosen so small 
that all critical points of polynomials  $g\in W$ 
are still non-degenerate, and that all critical points of $g$
corresponding to those critical points of $f$ converging to
periodic attractors, still converge  to periodic attractors.
(But it is possible that some additional critical points 
lie in basins of attracting periodic points of $g$.)
The other reason for the above construction (and for introducing
$\Sigma_0$) is that, provided we take  $r_0>0$ small, 
for each map $g\in \partial W_0$ there exists
a qc conjugacy defined on the basins $B(f)$ of periodic attractors
of $f$ to the basin $\tilde B(g)$
of the corresponding periodic attractor of $\tilde f$
with the additional property that the orbit of the critical points in 
$B(f)$ is mapped to the orbit of the corresponding critical points
in $\tilde B(g)$. (At this moment we don't say anything
about other periodic attractors that $g$ might have.)
Decreasing $W$ if necessary we can assume that 
all periodic points used in the construction 
of the complex box mapping $F\colon U\to V$
remain repelling for all polynomials from $W$
and that the conclusion Proposition~\ref{prop:nearbypol}
holds.

  Let $c_1,\ldots,c_N$ be the 
  critical points of $f$ whose iterates do not
  converge to periodic attractors (and so these points are in the
  Julia set of $f$) and which also are not mapped into the hyperbolic
  set mentioned above. As we have shown above, $N$ is the codimension
  of $QC(f)$ as a subspace of $\hat \Sigma(f)$ and is equal to the
  dimension of $W$.
   For $g\in W$ the critical point corresponding to
  $c_i$ will be denoted as $c_i(g)$. Suppose that for all $g\in W$ we
  have $g^k(c_i(g)) \neq c_1(g)$, $i=1,\ldots,N$, $k=1,\ldots$. Then
  all preimages of $c_1(g)$ move holomorphically with $g \in W$.
  Therefore we have a holomorphic motion on the set of all preimages
  of $c_1$ and, since this set is dense in the Julia set, we can
  extend the holomorphic motion to the whole Julia set. In fact, from   
  Theorem 7.4 of \cite{MS3} (which is based on 
  the harmonic $\lambda$-lemma),  any $g\in W$ is qc conjugate to $f$.
  But this contradicts the assumption that $QC(f)\cap \partial W_0=\emptyset$.

  So there exists $g \in W$ such that $g^k(c_i(g))=c_1(g)$ for some
  $k$ and $i$. If $N=1$ then $g$ is a semi-hyperbolic map and the theorem follows.
  Let $N$ be larger than one. The space $\Sigma_0$ is finite-dimensional and the
  equation $g^k(c_i(g))=c_1(g)$ determines an algebraic variety in
  this space. The singularities of this variety have complex
  codimension one in it and when we remove them we obtain a manifold
  of complex codimension one. Take a connected component of the
  intersection of this manifold with $W$, and denote this manifold
  by $M_{1}$.

  Now let us restrict ourselves to $M_{1}$ and consider all the
  preimages of $c_2(g)$, $g\in M_{1}$. Arguing as above we can
  deduce that either all maps in $M_1$ are qc conjugate or there is $g
  \in M_1$ and $i', k'$ such that $g^{k'}(c_{i'}(g))=c_2(g)$. In the
  second case, if $N=2$, we obtain a semi-hyperbolic map, 
  otherwise we can continue the construction. More precisely, we can consider a
  manifold $M_2$ which is obtained from the algebraic variety given by the
  equations $g^k(c_i(g))=c_1(g)$ and $g^{k'}(c_{i'}(g))=c_2(g)$
  removing singularities from it. The manifold $M_2$ has codimension
  two in $W$.

  Since the dimension of $W$ is equal to the number of critical points
  of $f$ which are in the Julia set ($=N$), we can continue until we
  obtain either a semi-hyperbolic map $g\in W$ or a manifold which
  consists of qc conjugate maps.
  In the first alternative, all remaining critical points of $g$
  either were already in the basin of a periodic attractor, or are
  mapped into a hyperbolic set. So as before, we again can perturb $g$
  to a hyperbolic map.  So assume the second alternative holds. Note
  that the manifold can consist of infinitely renormalizable maps, so
  we cannot use the Rigidity theorem directly. Also note that this
  manifold is obtained from an algebraic variety by removing
  singularities and intersecting this variety with $W$, so the
  intersection of the boundary of this manifold and the boundary of
  $W$ is not empty. This is because a complex algebraic variety cannot be
  contained in a bounded topological ball and if it has non empty
  intersection with such a ball, it also has non empty intersection
  with the ball's boundary.

  Obviously, we can construct infinitely many manifolds like this and,
  moreover, we can do it in such a way that these manifolds accumulate
  on $f$ (by taking a shrinking sequence of neighbourhoods $W_i\ni f$
  and finding maps $g$ as above in $W_i$). 
  Since these manifolds extend to the boundary of $W$ there
  exists $\tilde f\in \partial W_0$, $\tilde f\neq f$, such that some subsequence 
  of these manifolds
  accumulates on $\tilde f$ as well. So, we have obtained two sequence
  of maps $\{f_i\}$ and $\{\tilde f_i\}$ such that $f_i \to f$,
  $\tilde f_i \to \tilde f$ and such that for each $i$
  the maps $f_i$ and $\tilde f_i$ are qc conjugate.

For any $g\in W$ there exist 
complex box mapping $G\colon \hat U \to \hat V$
such that the domains of the box depend continuously
on $g$ as long as $g\in W$ (here we use Proposition~\ref{prop:nearbypol}). 
Applying this to $f_i,f, \tilde f_i, \tilde f$, and using the fact that
$\tilde F_i$ and $F_i$ are combinatorially equivalent, 
it follows that $F$ and $\tilde F$ are combinatorially equivalent.
  Hence, since $F$ is non-renormalizable, Theorem~\ref{thm:main2} 
   implies that $F$ and $\tilde F$ 
  are qc conjugate. It follows that $f$ and $\tilde f$ 
  are qc conjugate.
  Since $\tilde f\in \partial W_0$ we obtain
  a contradiction of the assumption that $QC(f)\cap \partial W_0=\emptyset$.

\section{Complex bounds implies rigidity for non-renormalizable 
complex box mappings}
\label{sec:rigboxbounds}

Assume that $f\colon U\to V$ and $\tilde f\colon \tilde U\to \tilde V$ are
 complex box mappings with only repelling periodic orbits,
 which are non-renormalizable
 and  combinatorially equivalent (in the sense defined
in Definition~\ref{def:equivbox}).
 
For simplicity, define $\mathcal{L}_x(P)$ to be the domain of the
first entry map to $P$ containing $x$ and let 
$\hat{\mathcal{L}}_{x}(P)$ is defined to be equal to $P$
if $x\in P$ and otherwise equal to $\mathcal{L}_x(P)$.

Given two critical puzzle pieces $P, Q$, we say that $Q$ is a 
{\it child} of $P$ if it is a unimodal pullback of $P$, i.e., 
if there exists a positive integer $n$ such that $f^{n-1}:f(Q)\to P$ is 
a diffeomorphism.
Given a puzzle piece $P\ni c$, by a {\it successor} of 
$P$ we mean a puzzle piece of the form $\hat{\mathcal{L}}_c(Q)$, 
where $Q$ is 
a child of $\hat{\mathcal{L}}_{c'}(P)$ for some $c'\in Crit(f)$.
The map $f$ is called {\em persistently recurrent}, if each critical 
puzzle piece has at most finitely many successors.

\begin{definition}
{\rm   A puzzle piece $P$ is called \emph{$\rho$--nice} if for any
  $x \in P \cap PC(f)$ one has $\mod(P - \LL_x(P)) \geq \rho$
  and \emph{$\delta$--fat} if there are puzzle
  pieces $P^+ \supset P \supset P^-$ such that the set 
  $P^+ - P^-$ does not
  contain points of the postcritical set of $f$, $\mod(P^+ - P)\geq\delta$
  and $\mod(P - P^-)\geq\delta$. 
  
  We say that a simply connected domain $U$ has 
  \emph{$\rho$--bounded geometry
    with respect to $x\in U$} if there are two disk $B(x,r)\subset U\subset B(x,R)$ and
  $R/r < \rho$.
  A domain $U$ is said to have \emph{$\rho$--bounded geometry} if there is $x\in U$
  such that $U$ has $\rho$--bounded geometry with respect to $x$.}
\end{definition}

\begin{theorem}[Strong complex bounds imply rigidity]\label{thm:rigid}
Assume that $f\colon U\to V$ and $\tilde f\colon \tilde U\to \tilde V$ 
are  non-renormalizable complex box mappings with all their periodic orbits repelling. Also assume that the map $H_0=H_1\colon V\setminus
U \to \tilde V\setminus \tilde U$ from
definition~\ref{def:equivbox} is qc.
Moreover, assume that the first return map $R_P$ 
to each critical puzzle piece $P$ of $f$ 
is either
\begin{description}
\item[a)\quad] non-persistently recurrent, or, 
\item[b)\quad] persistently recurrent and there exists a $\delta>0$ 
such that the complex 
box mapping $g:=R_P$ has arbitrarily small critical 
puzzle pieces $W$ which are $\delta$-nice, $\delta$-fat and
have $\delta$-bounded geometry with respect  to the critical point
in $W$; also assume that
the same statement holds 
for the corresponding puzzle pieces of $\tilde f$.
\end{description}
Then $f$ and $\tilde f$ are quasiconformally conjugate.
\end{theorem}
\begin{proof}
This was proved in Sections 6.3 and 6.4 of \cite{KSS1},
using  the QC-criterion from the Appendix of that paper
and the spreading principle from Section 5.3.
Note that since we assume here the conclusion rather than 
the assumption of the Key Lemma from \cite{KSS1},
the proof simplifies in a few places (for example in Lemma 6.7).
\end{proof}


In the persistently recurrent unicritical case, 
a different proof of the above theorem was recently 
given in \cite{AKLS} which does not require the
bounded geometry condition, but which does require
that the sequence of puzzle pieces $W_i$ from the second assumption
has the property that $W_{i+1}$ is a pullback of $W_i$
of uniformly bounded degree.

\bigskip

The next theorem asserts that the second assumption in the previous theorem
is always satisfied:

\begin{theorem}[Complex bounds for non-renormalizale complex box mappings]
\label{thm:boundsbox}
Assume that $f\colon U\to V$ is a  non-renormalizable complex box mapping with each of
its periodic orbits repelling. Moreover, assume that the first return map $R_P$ 
to each critical puzzle piece $P$ of $f$ 
is persistently recurrent.
Then there exists a $\delta>0$ 
such that the complex 
box mapping $g:=R_P$ has arbitrarily small critical 
puzzle pieces $W$ which are $\delta$-nice, $\delta$-fat and
have $\delta$-bounded geometry with respect to the critical point
in $W$.
The puzzle pieces for which this holds are combinatorially defined (the same statement holds 
for the corresponding puzzle pieces when $\tilde f$ is related to $f$
as in the previous theorem).
\end{theorem}

Theorem~\ref{thm:main2} now follows from the previous two theorems.
Indeed, if  the first return map to a critical puzzle piece
is  non-persistently recurrent, there exists a sequence
of critical puzzle pieces which under some iterates 
are mapped with bounded degree to a puzzle piece of fixed level. 
This implies in the non-persistently recurrent case 
that puzzle pieces shrink to points and absence of invariant linefields (and, in fact,
that the Julia set of the complex box mapping has zero Lebesgue measure).
In the persistently recurrent case, the previous theorem gives that
puzzle pieces shrink to points and absence of 
invariant linefields follows as in  
\cite[Theorem 10.3]{McM} and  \cite[Proposition 4.3]{Levin_Strien_loccon}.
In both cases, Theorem~\ref{thm:rigid} implies the qc rigidity of 
complex box mappings.

It follows that it suffices to prove Theorem~\ref{thm:boundsbox}. This will be done
in the remainder of this paper.

\section{Enhanced nest}
\label{sec:enest}
In the remainder of the paper we will prove
Theorem~\ref{thm:boundsbox}.  So let us assume that $f\colon U\to V$
is persistently recurrent.
In \cite{KSS1} we constructed a sequence of puzzle pieces 
\begin{equation}
\I_0\supset \I_1\supset \I_2 \ldots
\label{eq:enhan}
\end{equation}
around $c_0$, called the {\it enhanced nest} for the map $f$. 
The construction was based on the following lemma 
(Lemma 8.1 from \cite{KSS1}):

\begin{lemma}\label{lem:constrl}
  Let $\I\ni c$ be a puzzle piece. Then there is a positive integer $\nu$
  with $f^\nu(c)\in \I$ such that the following holds. Let
  $\U_0=\Comp_c(f^{-\nu}(\I))$ and $\U_j=f^j(U)$ for $0\leq j\leq \nu$. Then
  \begin{enumerate}
  \item $\#\{0\leq j\leq \nu-1: \U_j\cap\Crit(f)\not=\emptyset\}\leq b^2;$
  \item $\U_0\cap PC(f)\subset 
    \Comp_c\left(f^{-\nu}(\mathcal{L}_{f^{\nu}(c)}(\I))\right).$
  \end{enumerate}
\end{lemma}    

\medskip

For each puzzle piece $\I\ni c$, let $\nu=\nu(\I)$ be the smallest 
positive integer with the properties specified by Lemma \ref{lem:constrl}.
We define 
\begin{align*}
\mathcal{A}(\I)& =\Comp_c(f^{-\nu}(\mathcal{L}_{f^{\nu}(c)}(\I))),\\
\mathcal{B}(\I)& =\Comp_c(f^{-\nu}(\I)).
\end{align*}
As $f$ is persistently recurrent, each critical 
puzzle piece $\BP$ has a smallest successor, which we denote by
$\Gamma(\BP)$. Remark that if $\Q$ is an entry domain to $\BP$ intersecting 
$PC(f)$, then 
$\hat{\mathcal{L}}_c(\Q)$ is an successor  of $\BP$ by definition, 
and thus $\hat{\mathcal{L}}_c(\Q)\supset \Gamma(\BP)$.
Now define the enhanced nest (\ref{eq:enhan}) 
by  $\I_0=V_0$ and for each $n\geq 0$,  
$$\I_{n+1}=\Gamma^T\mathcal{B}\mathcal{A}(\I_n), 
$$
where $T=5b$. 
By construction $\I_{n+1}$ is a pullback of $\I_n$ and 
the map 
$$f^{p_n}\colon \I_{n+1}\to \I_n$$
has  degree bounded by  $M(d)$, where $M(d)$
only depends on the degree of $f$. 
This construction is chosen because of the following lemma
(see  Lemma 8.1 and 8.2 in \cite{KSS1}):

\begin{lemma} 
  \label{lem:transtime}
  For each $n>0$ there exists $\epsilon_n>0$ so that  
$\I_n$ is $\epsilon_n$-nice and $\epsilon_n$-fat. 
Moreover, denoting by  $r(\I_{n+1})$  
the minimal return time from $\I_{n+1}$ 
  to itself, we get
  \begin{enumerate}
  \item $3\,r(I_{n+1})\geq p_n$;   
  \item $p_{n+1}\geq 2p_n$.
    \end{enumerate}
\end{lemma}


Since the degree of $f^{p_n}\colon \I_{n+1}\to \I_n$ is bounded by the
$M(d)$, by Koebe's distortion lemma

\begin{lemma}
\label{lm:pullbackspace}
  There exists a universal constant $K$ (depending only on the degree $d$ of $f$) such
  that if $\I_n$ is $\epsilon$--nice then 
  $\I_{n+1}$ is $K\epsilon$--nice and $K\epsilon$--fat.
\end{lemma}

\section{Pullback lemmas}

Let $A$ be an annulus. We will use the two equivalent definitions of
its modulus.

Let $P_1(A)$ be the class of  non-negative Borel measurable function $\rho:A \to \R$
such that if $\gamma$ is any rectifiable Jordan closed curve separating boundaries
of $A$, then $\int_\gamma\rho \,d|z| > 1$. Then
$$ \mod(A) = \inf_{\rho\in P_1(A)} \int_A \rho^2 \,dz^2.$$
Similarly, let $P_2(A)$ be the class of 
non-negative Borel measure
function $\rho:A \to \R$ such that if $\gamma$ is any rectifiable Jordan curve
connecting boundaries of $A$, then $\int_\gamma\rho\,d|z| > 1$. Then
$$ \mod(A) = (\inf_{\rho\in P_2(A)} \int_A \rho^2 \,dz^2)^{-1}.$$
Write  $\C^+=\{z\in \C: Im(z) \geq0\}$.

\begin{lemma}[Small Distortion of Thin Annuli]
\label{lm:pullback}
  For every $K \in (0,1)$ there exists $\kappa>0$ such that if 
  $A\subset U$, $B\subset V$ 
  are simply connected domains symmetric with respect to the real
  line, 
  $F: U\to V$ 
  is a real holomorphic branched covering map of
  degree $D$ with all critical points real which can be decomposed as
  a composition of maps $F=f_1\circ\cdots \circ f_n$ with all maps $f_i$ 
  real and either
  real univalent or real branched covering maps with just one critical point,
  the domain 
  $A$ 
  is a connected component of 
  $f^{-1}(B)$ 
  symmetric
  with respect to the real line and the degree of 
  $F|_A$ 
  is $d$, then
  $$\mod(U-A) \geq \frac{K^D}{2d} \min\{\kappa, \mod(V-B)\}.$$
\end{lemma}

\begin{plm}
Let $B\subset V$ be two real-symmetric domains
and let $B^+= B\cap \C^+$. We claim that 
\begin{equation}
\frac 12 \mod(V-B^+) \leq \mod(V-B) \leq \mod(V-B^+).
\label{eq:halffactor}\end{equation}
The second inequality is obvious since $B^+\subset B$.
Let $\rho$ be in $P_1(V-B)$ such that $\int \rho^2 dz^2$ is almost $\mod(V-B)$.
Clearly, $\ \rho_0(z)=(\rho(z)+\rho(\bar z))/2$ is also in $P_1(V-B)$ 
and  $\int \rho_0^2  \,dz^2 \leq \int \rho^2 \,dz^2$
(where we use the inequality 
$((a+b)/2)^2\le (a^2+b^2)/2$). 
For $z\in V-B^+$ define $\rho_1$ by $\rho_1(z)=2\rho_0(z)$ if $Im(z)\geq0$
and $\rho_1=0$ otherwise. Obviously, $\rho_1\in P_1(V-B^+)$ and 
$$2\int \rho_0^2 \,dz^2 = \int \rho_1^2 \,dz^2 \geq \mod (V-B^+)$$
proving (\ref{eq:halffactor}).

Before we continue with the proof of Lemma~\ref{lm:pullback}
we state and prove:

\begin{sublemma}
There exists a universal constant $C>0$ such that the
following holds. Let $\D$ be a unit disk and $B\subset \D \cap \C^+$ be a simply connected
domain. Let $A$ be a connected component of $\phi^{-1}(B)$ where
$\phi(z)=z^d$.
Then
$$\mod(\D-A) > \frac{\mod(\D-B)}{1 + Cd\, \mod(\D-B)}.$$
\end{sublemma}
\begin{plm}
Let $\PP=\{z\in \C: \arg(z)\in [0,\pi/d]\}$ and
 $\QQ=\{z\in \C: \arg(z)\in [-\pi/(2d),3\pi/(2d)]\}$
(i.e. pullbacks by $z\mapsto z^d$ of the upperhalf plane
respectively the complex plane 
minus the negative imaginary axis).
Without loss of generality we can assume that $A$ is in the sector
$\PP$.  Let $\rho\in P_2(\D-B)$ and assume that  $\int \rho^2\, dz^2$
is close to $1 / \mod(\D-B)$. Pull $\rho$ back by $\phi$ to the domain $(\D-A)\cap \QQ$ and 
denote this pullback by $\rho_2$, and on $\D-\QQ$
 we set $\rho_2$ to be zero.  Define $\rho_3$ on $\D$
satisfying the following two properties: 
(i) for any Jordan rectifiable path
$\gamma$ connecting the set $[0,1/2 \exp(\pi i/d)] 
\cup [0, \frac 12]$ and the
boundary of $\D$ one has $\int_\gamma \rho_3\, d|z|>1$; 
(ii) for any Jordan
rectifiable 
curve $\gamma$ connecting the set $[\frac 12 \exp(\pi i/d),\exp(\pi i/d)]\cup [1/2,1]$  and $\C-\F$ one has
$\int_\gamma \rho_3\, d|z|>1$. 
It is easy to see that $\rho_3$ exists
with $C'=\int_\D \rho_3^2\, dz^2$ finite. 
If we take $\rho_3$  equal to $0$ outside the region
$F$ and constant inside, we see that
$C'/d\le C$ with $C$  independent of $d$.

Set $\rho_4(z)=\max\{\rho_2(z),\rho_3(z)\}$ for $z\in \D-A$. It is easy to check
that $\rho_4\in P_2(\D-A)$. Then
$$\frac 1{\mod(\D-A)} \leq \int \max\{\rho_2(z)^2,\rho_3(z)^2\} \, dz^2 \leq 
$$
$$\quad \quad \leq \int
\rho_2(z)^2 \, dz^2 + \int \rho_3(z)^2 \, dz^2\leq \int \rho_2(z)^2 \, dz^2 +C'.$$
This implies
$$\mod(\D-A)\geq\frac{\left[\int \rho_2(z)^2 \, dz^2\right]^{-1}}
{1+C'\left[\int \rho_2(z)^2 \, dz^2\right]^{-1}} \approx \frac
{\mod(\D-B)}{1+C'\mod(\D-B)}.$$
\end{plm}

Lemma~\ref{lm:pullback} is a direct consequence of 
equation (\ref{eq:halffactor})
and the sublemma. 
First, we split $B$ into two parts and then only look at $B^+$.
Because of the second inequality in (\ref{eq:halffactor}), we do not lose any modulus at this point. Then we pullback $B^+$ 
many times all the way. 
At each step we do not lose modulus if the pullback is univalent. 
If the $i$-th pullback is quadratic but the corresponding pullback 
of $B^+$ does not contain a critical value of the map $f_i$, we apply
the sublemma and  so lose just a bit of
modulus of $\D-B^+$ which is small; in this case set $e_i=1$. 
If the corresponding 
pullback of $B^+$ does contain a critical value of $f_i$ of order 
$d_i$, then we lose
a factor $e_i=d_i$.  At the end we reconstruct the preimage 
of $B$ by mirroring the preimage of $U^+$ and use the 1st inequality  in (\ref{eq:halffactor}). This means that we lose 
a factor $2$. So in total we lose a factor $K^D/(2d)$, since
$d=e_1\dots e_n$ where $e_i=1$ or $e_i=d_i$.  
\end{plm}

In the real case the previous lemma implies and sharpens the
following result of J.Kahn and M.Lyubich, see \cite{KL}:

\begin{lemma}\label{lm:KL}
  For any $\eta > 0$ and $D>0$ there is $\epsilon=\epsilon(\eta,D)>0$ such that the
  following holds: Let $A\subset A' \subset U$ and $B\subset B'\subset V$ be topological disks
  in $\C$ and let $F: (A, A', U)  \to (B, B', V)$ be a holomorphic
  branched covering map. Let the degree of $f$ be bounded by $D$ and
  the degree of $f|_{A'}$ be bounded by $d$. Then
  $$\mod(U\setminus A) > \min(\epsilon,\, \eta^{-1} \mod(B'\setminus B),\, C\eta d^{-2} \mod(V\setminus
  B)),$$
  where $C>0$ is some universal constant.
\end{lemma}

\section{Complex bounds for the enhanced nest}
\label{sec:complexboundsenhanced}

In this section we will prove the following proposition

\begin{proposition}[Complex bounds]\label{thm:complexbounds}
  Let $f$ be a complex box mapping.
  For any $\epsilon>0$ there exists $\delta>0$ such that if $\I_0$ is $\epsilon$--nice,
  then all $\I_n$, $n=1,2,\ldots$, are $\delta$--nice.
\end{proposition}

Because of Lemma~\ref{lm:KL}, $\I_n$, $n=2,3,\dots$ 
are also $\delta'>0$-fat.

\medskip

\begin{prf} Denote $\mu_n$ to be such that $\I_n$ is $\mu_n$-nice.
Fix some integer $M>4$ and suppose that $n>M+1$. Let 
 $A\subset\I_n$ be some
  domain of the first return map to $\I_n$ containing a point of the
  postcritical set and let $r$ be its return time.

  {\bf Step 1:} $A$ is contained in $\LL_x(\I_{n-4})$ for some $x$.
  Indeed, write   $$P_{n,M}={p_{n-1}+p_{n-2}+\ldots+p_{n-M}}
  \mbox{ and }B:=f^{P_{n,M}}(A).$$
  Note that $f^{P_{n,M}}$ maps $\I_n$ onto $\I_{n-M}$. 
  By the last two inequalities in 
  Lemma~\ref{lem:transtime},
  $$r\ge r(\I_n)\ge \frac{1}{3}p_{n-1}\ge
  \frac{16}{3}p_{n-5} \ge p_{n-5}+p_{n-6}+p_{n-7}+\dots
  \ge P_{n-4,M-4}$$ 
  and so $s:=r-p_{n-4,M-4}\ge 0$. Since
  $\I_n=f^r(A)=f^s\circ f^{n-4,M-4}(A)$ we get that
  $f^s(f^{p_{n,M}}(A))$
  is equal to 
  $$f^{p_{n-1}}\circ \dots 
  \circ f^{p_{n-4}}\circ f^s\circ f^{n-4,M-4}(A)=\I_{n-4}
  $$
  and therefore $f^{P_{n,M}}(A)$ is contained in some
  $\LL_x(\I_{n-4})$.

  {\bf Step 2:} There exists $\I_k^+$ with 
  $\I_k\subset\I_k^+\subset \I_{k-1}$ such that
  $\mod(\I_k^+ - \I_k)> K_1\mu_{k-1}$ and $\I_k^+ -\I_k$ does not
  contain postcritical points of $f$. 
 Indeed,  the domain $\I_k$ is $K_1\mu_{k-1}$--fat, where $K_1$ is a constant given by Lemma~\ref{lm:pullbackspace}.
  Since $\I_{k-1}$ is $\mu_{k-1}$
  nice as well we can assume that $\I_k^+\subset\I_{k-1}$.
  
  {\bf Step 3:} 
  Fix $x\in f^{P_{n,M}}(A) \cap PC(f)$ 
  and some integer $k \in [n-4, n-M]$.
  Take $\nu$ so that
  $R_{\II_{k}}|_x=f^{\nu}$. Denote the pullback of $\I_k^+ - \I_k$ by
  $f^\nu$ by $\A_k$, \ie the point $x$ is surrounded by the annulus $\A_k$
  and $f^\nu(\A_k)=\I_k^+ - \I_k$.  The degree $d$ of
  $f^\nu|_{\LL_x(\I_k)}$ is bounded by some universal constant $K_2$
  depending only on $b$ and the degree of the map
  $f^\nu|_{\Comp_x(f^{-\nu}(\I_k^+))}$ is also $d$. Hence, 
  \begin{equation}\mod (\A_k) > (K_1/K_2) \mu_{k-1}.
  \label{eq:akineq}
  \end{equation}

 {\bf Step 4:} Obviously, all annuli $\A_k$, $k=n-4,\ldots,n-M$, are nested and surround
  $f^{P_{n,M}}(A)$. This implies that 
 \begin{equation} \mod (\I_{n-M}-f^{P_{n,M}}(A)) > (K_1/K_2) (\mu_{n-M-1}+\ldots+\mu_{n-5}).
  \label{eq:sum}\end{equation}

 {\bf Step 5:} The degree of the map $f^{P_{n,M}}|_A$ is bounded by some constant
  $d$ which does not depend on $M$, while the degree of
  $f^{P_{n,M}}|_{\I_n}$ is bounded by some constant $K_4(M)$ depending on $M$.
  The second assertion is obvious, so let us prove the first one.
  Indeed, decompose the map $f^{P_{n,M}}|_A$ as
  $f^{P_{n-8,M-8}}\circ f^{P_{n,8}}$. The degree of $f^{P_{n,8}}$ is
  bounded by some constant depending only on $b$. The domain
  $f^{P_{n,M}}(A)$ is contained in $\LL_x(\I_{n-4})$ and $M>8$, hence
  $f^{P_{n,8}}(A)$ belongs to some other component
  $\LL_y(\I_{n-4})$. Due to Lemma~\ref{lem:transtime} we know that
  $r_{n-4}>P_{n-8,M-8}$, therefore the degree of the map
  $f^{P_{n-8,M-8}}$ is less or equal the degree of $R_{\I_{n-4}}|_{\LL_y(\I_{n-4})}$
  which is bounded by some constant depending only on $b$.
  
 {\bf Step 6:}    Now apply Lemma~\ref{lm:KL}. Let $B$ and $B'$ be
  bounded by 
  the inner and outer boundary of the annulus $A_{n-4}$ and let $V:=
  \I_{n-M}$. Let $A$, $A'$, $U$ be the corresponding pullbacks of $B$,
  $B'$, and $V$ by $F^{-1}$ where $F:=f^{P_{n,M}}$. 
  Notice that the degree of maps $f^{P_{n,M}}|_A$ and
  $f^{P_{n,M}}|_{A'}$ are the same because $A_{n-4}$ does not contain
  postcritical points.
  

  We may assume that $K_1/K_2<1$.
  Fix $\eta=\frac{K_1}{2K_2}$ and let $C$ and $\epsilon(\eta,D)$ 
  be the 
  the constants from  Lemma~\ref{lm:KL} associated 
  to $\eta$ and $D$. Now fix $M=\frac 4C d^2 K_2/K_1 +6$
  and take $d$ 
  and $D=K_4(M)$ as in Step 5.
  From equation~(\ref{eq:sum}),
  \begin{equation}
  \mod(V\setminus B) > \frac{K_1}{K_2} (\mu_{n-M-1}+\cdots +\mu_{n-5})
  \ge \frac{4}{C} d^2\mu_{n-5,M-5}
  \label{eq:vbineq}
  \end{equation}
  where $\mu_{k,m}=\min(\mu_k,\mu_{k-1},\ldots,\mu_{k-m})$.
   Lemma~\ref{lm:KL} together with 
   equations (\ref{eq:akineq}) and (\ref{eq:vbineq}) 
   gives
  $$\mod(U \setminus A) > \min\left(\epsilon(\eta,M),\, \eta^{-1}\mod(B'\setminus B),\, 
  C\eta d^{-2}\mod(V\setminus B) \right)$$
  $$>  \min\left(\epsilon(\eta,M),\, 2\mu_{n-5},\, 2\mu_{n-5,M-5}\right),$$
  and so $\mu_n> \min(\epsilon(\eta,M), 2\mu_{n-5,M-5})$. By
  Lemma~\ref{lm:pullbackspace} we also have $\mu_n > K\mu_{n-1}$. These
  two inequalities prove the proposition.
\end{prf}

If $f$ is real we can apply Lemma~\ref{lm:pullback} instead of
 Lemma~\ref{lm:KL}. This makes the above proof slightly easier.

\section{Bounded Geometry (Proof of Theorem~\ref{thm:boundsbox})}
\label{sec:bounded-geometry}

In this section we will show that the above complex 
bounds immediately imply 
that the puzzle pieces from the enhanced nest 
have bounded geometry.

\begin{proposition}[Bounded geometry]\label{keylemma}
Let $f$ be a complex box mapping (not necessarily real).
Let $\I_0$ be $\epsilon$--nice and let it has $\rho$--bounded geometry with
  respect to $c_0$.
 Then all $\I_n$, $n=1,2,\ldots$, have $\delta$--bounded geometry, where $\delta$
 depends only on $\rho$ and $\epsilon$.
 \end{proposition}

Let us first state and prove the following easy consequences of 
 Koebe's distortion lemma:

\begin{lemma}
\label{lm:bg1}
  Let a domain $U$ have $\rho$--bounded geometry with respect to some
  point $x$ and let $A\subset U$ be a domain containing $x$. Then 
  $U$ has
  $K\rho$--bounded geometry with respect to all $y\in A$, where the constant
  $K$ depends only on $\mod(U-A)$.
\end{lemma}
\begin{plm}
By Koebe's distortion lemma it follows that $d(x,y)\le K d(y,\partial U)$ for any $y\in A$
where $K$ depends on $\mod(U-A)$.
So if $d(x,y)\ge (1/2) d(x,\partial U)$ then
$$\frac{\sup_{z\in \partial U}d(y,z)}{d(y,\partial U)}
\le 2K \frac{\sup_{z\in \partial U}d(x,z)}{d(x,y)}
\le 4K \frac{\sup_{z\in \partial U}d(x,z)}{d(x,\partial U)}$$
and we are done.
If $\diam(A)\le (1/2) d(x,\partial U)$ then
$d(y,\partial U)\ge (1/2)d(x,\partial U)$ and the assertion
follows immediately (without using Koebe's distortion lemma).
\end{plm}

Similarly we have

\begin{lemma}
\label{lm:bg2}
  Let $f:U \to V$ be a holomorphic covering map, $B\subset V$, 
  $A$ is a connected component of $f^{-1}(B)$.  Then if 
  $B$ has $\rho$--bounded
  geometry with respect to some point $y\in B$, then $A$ has 
  $K\rho$--bounded geometry with respect to $x$, where the point 
    $x\in A$ is any
  preimage of $y$ by $f$ and the constant $K$ depends only on
  $\mod(V-B)$ and on the degree of the map $f$.
\end{lemma}

\noindent
\begin{proofof}{Proposition~\ref{keylemma}}
Let $\I_n$ have $\rho_n$--bounded geometry with respect to $c_0$. Since
we know that $\I_n$ is $\delta$--fat for some constant $\delta$, $\I_n$
has $K_1\rho_n$--bounded geometry with respect to $f^{p_n}(c_0)$ where
$K_1$ is given by Lemma~\ref{lm:bg1}. Then $f(\I_{n+1})$ has
$K_1K_2\rho_n$--bounded geometry with respect to $f(c_0)$ where $K_2$ is
given by Lemma~\ref{lm:bg2} (notice that the degree of the map
$f^{p_n-1}$ depends only on $b$). Therefore, $\I_{n+1}$ has
$\sqrt{K_1K_2\rho_n}$--bounded geometry with respect to $c_0$, \ie
$\rho_{n+1}\geq\sqrt{K_1K_2\rho_n}$.
\end{proofof}

Because of Propositions~\ref{thm:complexbounds} and \ref{keylemma}, we have completed the proof 
of Theorem~\ref{thm:boundsbox}.

\bibliographystyle{plain}
\bibliography{bibnewrigid.bib}


\end{document}